\documentclass[12pt]{amsart}
\usepackage{amssymb, mathrsfs,amsmath}
\usepackage{latexsym}
\usepackage{amsfonts}
\usepackage{shadow}
\date{}
\author[D. Alpay]{Daniel Alpay}
\author[H. Attia]{Haim Attia}
\address{(DA, HA) Department of Mathematics \newline
Ben Gurion University of the Negev \newline P.O.B. 653, \newline
Be'er Sheva 84105, \newline ISRAEL} \email{dany@math.bgu.ac.il,
atyah@bgu.ac.il}
\author[D. Levanony]{David Levanony}
\address{(DL) Department of Electrical Engineering \newline
Ben Gurion University of the Negev \newline P.O.B. 653, \newline
Be'er Sheva 84105, \newline ISRAEL }
\thanks{D. Alpay thanks the
Earl Katz family for endowing the chair which supported his
research. The research of the authors was supported in part by the
Israel Science Foundation grant 1023/07}
\email{levanony@ee.bgu.ac.il} 
\title[White noise based stochastic calculus]
{White noise based stochastic calculus associated
with a class of Gaussian processes}
\begin{document}
%\maketitle
\parindent 0cm
\newtheorem{Pa}{Paper}[section]
\newtheorem{Tm}[Pa]{{\bf Theorem}}
\newtheorem{La}[Pa]{{\bf Lemma}}
\newtheorem{Cy}[Pa]{{\bf Corollary}}
\newtheorem{Rk}[Pa]{{\bf Remark}}
\newtheorem{Pn}[Pa]{{\bf Proposition}}
%\newthePb}[Pa]{{\bf Problem}}
\newtheorem{Dn}[Pa]{{\bf Definition}}
\newtheorem{Ex}[Pa]{{\bf Example}}
\numberwithin{equation}{section}
\def\L{\mathbf L}
\def\R{\mathbb R}
\def\N{\mathbb N}
\def\C{\mathbb C}
\def\s{\mathscr S}
\def\ss{\mathscr S^\prime}
\def\sr{\mathscr S(\R)}
\def\ssr{\mathscr S'(\R)}
\def\(s){\mathscr S(\R^2)}
\def\F{\mathcal F}
\def\P{\mathcal P}
\def\W{\mathcal W}
\def\Dom{{\rm dom}~(T_m)}
\def\Doms{{\rm dom}~(T_m^*)}
\def\Def{\stackrel{{\rm def.}}{=}}
\begin{abstract}
Using the white noise space setting, we define and study
stochastic integrals with respect to a class of stationary
increment Gaussian processes. We focus  mainly on continuous
functions with values in the Kondratiev space of stochastic
distributions, where use is made of the topology of nuclear
spaces. We also prove an associated Ito formula.
\end{abstract}
\maketitle
\tableofcontents
\section{Introduction}
\setcounter{equation}{0} In this paper we study stochastic
integration with respect to stationary increment Gaussian
processes $\left\{X_m(t)\,,\, t\in\mathbb R\right\}$ with
covariance functions of the form
\begin{equation}
\label{cm}
\begin{split}
C_m(t,s)\stackrel{\rm def.}{=}
E[X_m(t)\overline{X_m}(s)]&=\int_{\mathbb
R}\frac{(e^{iut}-1)(e^{-ius}-1)}{u^2}m(u)du\\
&=r(t)+\overline{r(s)}-r(t-s)-r(0),
\end{split}
\end{equation}
where $m$ is a measurable positive function subject to
\begin{equation} m(u)\leq
\begin{cases}
K\left|u\right|^{-b}\quad{\rm if}\quad|u|\leq 1,\\
K|u|^{2N}\hspace{1mm}\quad{\rm if}\quad|u|>1,
\end{cases}
\label{mbound}
\end{equation}
with $b<2$, $N\in{\mathbb N}_0$, and $0<K<\infty$, where
\[
r(t)=-\int_{\R} \Big\{e^{itu}-1-\frac{itu}{u^2+1}
\Big\}\frac{m(u)}{u^2}du.
\]
An associated
Ito formula is subsequently derived.
\bigskip

We use the white noise space setting as developed by T. Hida, and
in particular the Gelfand triple $(S_1,\mathcal W, S_{-1})$
consisting of the Kondratiev space $S_1$ of stochastic test
functions, of the white noise space $\W$,  and  the Kondratiev
space $S_{-1}$ of stochastic distributions; see \cite{MR1244577},
\cite{MR2444857}, \cite{MR1408433}. Various notions pertaining to
these works, which are used in the
introduction, are recalled in Sections \ref{prelim} and \ref{wns}.\\

Explicit constructions of $X_m$ and its derivative, which we use
below, are detailed in \cite{aal2} utilizing this setting. A key
role in the arguments of \cite{aal2} is played by the operator
\begin{equation}
\widehat{T_mf}(u) {\stackrel{{\rm
def.}}{=}}\sqrt{m(u)}\widehat{f}(u), \label{Tm}
\end{equation}
where $\widehat{f}$ denotes the Fourier transform of $f$:
\[
\widehat{f}(u)=\int_{\mathbb R}e^{-iux}f(x)dx.
\]

Since $m$ satisfies \eqref{mbound}, the domain of $T_m$
%\begin{equation}
%\label{voltaire} {\rm dom}~(T_m){\stackrel{{\rm def.}}{=}}
%\left\{f\in {\mathbf L}_2(\R):~~\int_{\R}m(u) \left|\widehat{f}(u)
%\right|^2du<\infty\right\},
%\end{equation}
contains in particular the Schwartz space $\sr$.  We note that the operator $T_m$ will
not be local in general: The support of $T_mf$ need not be
included in the support of $f$. The example
\begin{equation}
\label{u4}
m(u)=u^4e^{-2u^2},
\end{equation}
given in \cite{aal2} illustrates this point. The choice
\begin{equation}
\label{BH}
m(u)=\frac{1}{2\pi}|u|^{1-2H}du,\quad H\in(0,1),
\end{equation}
corresponds to the fractional Brownian motion $B_H$ with Hurst
parameter $H$, such that
\[
E(B_H(t)B_H(s))=V_H\left\{|t|^{2H}+|s|^{2H}-|t-s|^{2H}\right\},
\]
where
\begin{equation}
\label{VH}
V_H= \frac{\Gamma(2-2H)\cos(\pi H)}{\pi(1-2H)H},
\end{equation}
with $\Gamma$ denoting the Gamma function. For this choice of $m$,
the operator $T_m$ has been introduced in \cite[(2.10),
p. 304]{eh} and in
\cite[ Definition 3.1, p. 354]{bosw}.\\

It is easy to see that, when $m$ is even,
\begin{equation}
\overline{T_mf}=T_m\overline{f}.
\end{equation}
In this paper we focus on the real-valued case, and therefore will restrict
ourselves to even functions $m$.\\

Let $Y$ be an $S_{-1}$-valued continuous function defined for
$t\in[a,b]$. Our main result, see Theorem \ref{bastille} below,
states that the integral
\begin{equation}
\label{sto_int}
\int_{\mathbb R}Y(t)\lozenge W_m(t)dt,
\end{equation}
(with $\lozenge$ being the Wick product to be defined below)
usually understood in the sense of Pettis,
%where the second integral is understood in the sense of Pettis. When $m$
%is equal to \eqref{BH}, this reduces to the definition given in
%\cite{bosw}, \cite{}.\\
is a limit of Riemann sums, with convergence in a Hilbert space
norm.

% when $Y$ is continuous and defined on a compact

%interval.

In the case of the fractional Brownian motion, a related
characterization was given in \cite[(3.16), p. 591]{MR1741154}.
Still, our methods
and the methods of \cite{MR1741154} are quite different.\\

%We also consider the case where $Y$ is adapted with respect to

%the white noise and  prove a Ito-type formula for $X_m$.\\

The paper consists of six sections besides the introduction.
With the proof of Theorem \ref{bastille} in mind, we begin
Section \ref{prelim} with a short review of the topology of
countably normed spaces and of their duals. In the next section,
we prove a new result on continuous functions with valued in the
dual of a perfect space. The main features of Hida's theory of
the white noise space are then reviewed in Section \ref{wns}.
Various notions, such as the Wick product and the Kondratiev
spaces, appearing in this introduction, are defined there.
In Section \ref{wickito} we state and prove
Theorem \ref{bastille}. An Ito-type formula is proved in Section
\ref{itoform}. The last section is devoted to a number of
concluding observations.\\

%The connection with the Wick-Skorokhod integral is done in

%Section \ref{skho}.

Notation is standard. In particular we set
\[
{\mathbb N}=\left\{1,2,3,\ldots\right\} \quad{\rm and}\quad
{\mathbb N}_0={\mathbb N}\cup\left\{0\right\}.
\]

\section{Countably normed spaces}
\setcounter{equation}{0} \label{prelim} Nuclear spaces are an
indispensable part of the foundation upon which white noise
theory, to be utilized below, is built.  In this section we
review part of the theory of nuclear spaces, as developed in
\cite{GS2} and \cite{MR35:7123}. We use the notation of these
books.\\

Let $\Phi$ be a vector space (on $\mathbb R$ or $\mathbb C$)
endowed with a sequence of norms $(\|\cdot\|_p)_{p\in\mathbb N}$.
Assume that the norms are defined by inner products and that the
sequence is increasing:
\[
p\le q\Longrightarrow \|h\|_p\le \|h\|_q,\quad \forall h\in\Phi.
\]

Denote by $\mathcal H_p$ the closure of $\Phi$ with respect to the
norm $\|\cdot\|_p$. For $p\le q$, a Cauchy sequence in $\mathcal
H_q$ is a Cauchy sequence in $\mathcal H_p$, and this defines a
natural map from $\mathcal H_q$ into $\mathcal H_p$. In general,
this map need not be one-to-one. A counterexample is presented in
\cite[p. 13]{GS2}. This phenomenon will not occur in the case of
reproducing kernel Hilbert spaces, as is shown in the following
proposition. In this statement, recall that the positive kernel
$K_q$ is said to be smaller than the positive kernel $K_p$ if the
difference $K_p-K_q$ is positive.

\begin{Pn}
\label{kernel}
Given the notation above, let $p\le q$, and
assume that $\mathcal H_p$ and $\mathcal H_q$
are reproducing kernel Hilbert spaces of functions defined on a common
set $\Omega$, with respective reproducing kernels $K_p$ and $K_q$.
Assume that
\[
K_q(z,w)\le K_p(z,w)
\]
in the sense of reproducing kernels. Then, $\mathcal H_q$ is a subset of
$\mathcal H_p$, and the inclusion is contractive.
\end{Pn}

{\bf Proof:} This follows from the decomposition
\[
K_p(z,w)=K_q(z,w)+(K_p(z,w)-K_q(z,w))
\]
of $K_p$ into a sum of two positive kernels, and of the
characterization of the reproducing kernel Hilbert space
associated with a sum of positive kernels. See \cite[\S 6]{aron}
for the latter.
\mbox{}\qed\mbox{}\\

In the sequel we assume that $\mathcal H_q\subset \mathcal H_p$
when $p\le q$. The inclusion will not be in general an isometry.
The space $\Phi$ is the projective limit of the spaces $\mathcal
H_p$. It will be complete if and only if
\[
\Phi=\bigcap_{n=1}^\infty\mathcal H_n.
\]
See \cite[Th\'eor\`eme 1, p. 17]{GS2}.\\

Denote by $\Phi^\prime$ the topological dual of $\Phi$. Then
\[
\Phi^\prime=\bigcup_{n=1}^\infty{ \mathcal H}_n^\prime,
\]
where ${\mathcal H}_n^\prime$ denotes  the topological dual of
${\mathcal H}_n$. Furthermore, denote by
\[
\langle v,u\rangle, \quad v\in\Phi^\prime,\quad u\in\Phi,
\]
the duality between $\Phi$ and $\Phi^\prime$. By definition, for
$u\in\mathcal H_r$ and $v\in\mathcal H_r^\prime$ one has
\[
\langle v,u\rangle=\langle v,u\rangle_r,
\]
where $\langle v,u\rangle_r$ denotes the duality between
$\mathcal H_r$ and $\mathcal H_r^\prime$, and
\begin{equation}
\label{eq:ieqrtyu} \|v\|_{{\mathcal
H}^\prime_r}=\sup_{\substack{u\in{{\mathcal H}_r},\\
\|u\|_{{\mathcal H}_r}=1}}\langle v,u\rangle_r, \quad{\rm
and}\quad
%\label{eq:ineqrty}
|\langle v,u\rangle_r|\le\|v\|_{{\mathcal H}^\prime_r}\|u\|_{{\mathcal H}_r}.
\end{equation}
Moreover, note that, for $p\ge r$ and $v\in\mathcal H_p^\prime$
and $u\in\mathcal H_r\subset\mathcal H_p$, one has:
\begin{equation}
\label{inner_pro_equal}
\langle v,u\rangle=\langle v,u\rangle_r=\langle v,u\rangle_p.
\end{equation}

%\begin{Dn}

%\label{prelim1}

%.

%\end{Dn}

Indeed, \eqref{inner_pro_equal} expresses the valued of the
linear functional $v$ on $u$.
See \cite[p. 56]{MR35:7123} for a discussion of this point.\\

We refer the reader to \cite[\S 5.1, pp. 41-44]{GS2} for the definition of
the strong topology on $\Phi^\prime$. To ease the reading of
Gelfand-Shilov \cite{GS2}, we make the following remark: In
verifying that a topological vector space $V$ is Hausdorff, it is
necessary and sufficient to check the following: For every $v\in
V$ there exists a neighborhood of $0$, say $\mathcal N$, such that
$v\not\in \mathcal N$. See for instance \cite[Proposition 9, p.
70]{MR0162111}. This is the condition which is used and
verified in \cite[\S 5.1]{GS2}.\\

The complete, countably normed space $\Phi$ is called {\sl
perfect}, or {\sl a Montel space}, when any subset of $\Phi$ is
bounded and closed if and only if it is compact. A necessary
condition for $\Phi$ to be perfect is that, for every
$r\in\mathbb N$ there exists a $p>r$ such that the inclusion from
$\mathcal H_p$ into $\mathcal H_r$ is compact. See
\cite[Th\'eor\`eme 1, p. 55]{GS2}. It will be called  nuclear if
the above inclusion can be chosen to be of trace class (as an
operator between Hilbert spaces).

\section{Continuous functions
with values in the dual of a perfect space}

\setcounter{equation}{0}
The following theorem is the key in our construction of the
stochastic integral as a limit of Riemann sums.

\begin{Tm}
\label{function}
Let $(E,d)$ be a compact metric space, and let
$f$ be a continuous function from $E$ into the dual of a
countably normed perfect space $\Phi=\cap_{n=1}^\infty \mathcal
H_n$, endowed with the strong topology. Then there exists a
$p\in\mathbb N$ such that $f(E)\subset \mathcal H_p^\prime$, and
$f$ is uniformly continuous from $E$ into $\mathcal H_p^\prime$,
the latter being endowed with its norm induced topology.
\end{Tm}

{\bf Proof:} We divide the proof into a number of steps.\\

STEP 1: {\sl $f(E)$ is compact.}\\

Indeed,  the dual  space $\Phi^\prime$ endowed with the
strong topology is a Hausdorff space; see  \cite[\S 5.1, pp. 41-42]{GS2}.
Since $E$ is compact and the function $f$ is continuous,
it follows that $f(E)\subset \Phi^\prime$ is compact. \\

STEP 2: {\sl  There exists an $r\in\mathbb N$ such that
$f(E)\subset\mathcal{H}_{r}^\prime$ and $f(E)$ is bounded in
$\mathcal{H}_{r}^\prime$.}\\

Since $f(E)$ is compact, it is bounded in $\Phi^\prime$; See
\cite[Proposition 1, p. 54]{GS2}, and, see \cite[D\'efinition 5,
p. 30]{GS2} for the notion of a bounded set in a topological
vector space. By the characterization of bounded sets in the
strong dual of a perfect space,  see \cite[Th\'eor\`eme 2 p.
45]{GS2}, there exists an $r\in\N$ such that
$f(t)\in\mathcal{H}_{r}^\prime$ for all $t\in E$.\\

%This follows from \eqref{eq:ieqrtyu}.\\

STEP 3: {\sl  Set $r$ as in the previous step. Let $t\in E$ and
let $(s_m)_{m\in\mathbb N}$ be a sequence of elements of $E$ such
that $\lim_{m\rightarrow\infty}d(t,s_m)=0$. Then, for every
$h\in\mathcal H_r$,}
\begin{equation}
\label{eq:chaos}
\lim_{m\rightarrow\infty}\langle f(s_m)-f(t),h\rangle=
\lim_{m\rightarrow\infty}\langle f(s_m)-f(t),h\rangle_r=0.
\end{equation}

Indeed, the function $f$ is continuous in the strong topology of
$\Phi^\prime$, and therefore sequentially
continuous in the strong topology, and hence weakly sequentially continuous. \\

STEP 4: {\sl There exists  a $p>r$ such that the inclusion map
from
$\mathcal H_p$ into $\mathcal H_r$ is compact.}\\

Such a $p$ exists since the space $\Phi$ is  assumed perfect. \\

Before turning to the fifth step, we remark the following: Since
$E$ is a metric space, it is enough to verify continuity by using
sequences; see for instance \cite[Th\'eor\`eme 4 p.
58]{MR0162111}. Furthermore, we note that in Step 3, we cannot say
in general that $f(s_n)$ tends to $f(t)$ in the

norm of $\mathcal H_r^\prime$, yet we have:\\

%Furthermore

%\[

%\lim_{n\rightarrow\infty}\langle f(t_{n_m})-f(t),h\rangle=0

%\]

%for all $h\in\mathcal H_r$.\\

STEP 5: {\sl   Set $p$ as in Step 3. The function $f$ is
continuous from $E$ into
$\mathcal H_p^\prime$, the latter endowed with its
norm topology.}\\

Set $t\in E$ and let $(t_n)_{n\in\mathbb N}$ be a sequence of
elements of $E$ such that $\lim_{n\rightarrow\infty}d(t_n, t)=0$.
Since $f$ is continuous, it follows
\begin{equation}
\label{eq:cont}
f(t_n)\to f(t)
\end{equation}
in the strong topology of $\Phi^\prime$. Using
\cite[Th\'eor\`eme 4, p. 58]{GS2}, one has $ f(t_n)\to f(t) $ in
norm in $\mathcal{H}_{p}^\prime$. In the proof of
\cite[Th\'eor\`eme 4, p. 58]{GS2}, the integers $r$ and $p$ depend
{\it a priori} on the sequence. We repeat this argument and verify
that the {\it same} $r$ and $p$ can be taken for all sequences $f(t_n)$:\\

The argument of \cite{GS2} is as follows.  Assume that
\eqref{eq:cont} does not hold. Then, there exist an $\epsilon>0$,
a subsequence $(t_{n_m})_{m\in\mathbb N}$, and a sequence
$(h_m)_{m\in\mathbb N}$ of elements in the closed unit ball of
$\mathcal H_p$  such that
\begin{equation}
\label{ep}
|\langle f(t_{n_m})-f(t),h_m\rangle_p|\ge \epsilon.
\end{equation}

Since the inclusion map from $\mathcal H_p$ into $\mathcal H_r$
is compact, the sequence $(h_m)_{m\in\mathbb N}$ has a convergent
subsequence in $\mathcal H_r$. Denote this subsequence by
$(h_m)_{m\in\mathbb N}$ as well, and set
\[
h=\lim_{m\rightarrow\infty} h_m\in\mathcal H_r.
\]

Writing (recall \eqref{inner_pro_equal})
\[
\langle f(t_{n_m})-f(t),h_m\rangle_p=\langle f(t_{n_m})-f(t),h_m-h\rangle_r+
\langle f(t_{n_m})-f(t),h\rangle_r ,
\]
we see that
\[
\lim_{m\rightarrow\infty}  \langle f(t_{n_m})-f(t),h_m\rangle_p=0.
\]

Indeed, using \eqref{eq:ieqrtyu}, we have
\[
\lim_{m\rightarrow\infty}
\langle f(t_{n_m})-f(t),h_m-h\rangle_r=0,
\]
since $f(E)$ is bounded in $\mathcal H_r^\prime$, and from Step 3,
with $s_m=t_{n_m}$,
\[
\lim_{n\rightarrow\infty}\langle f(t_{n_m})-f(t),h\rangle=0.
\]

A contradiction with \eqref{ep} is thus obtained, thus verifying
the STEP 5 statement.\\

STEP 6: {\sl $f$ is uniformly continuous from $E$ into $\mathcal
H_p^\prime$.}\\

This stems from the fact that $E$ is compact and that
 $\mathcal H_p^\prime$ is Hausdorff.\\

%

% Voir si necessaire Garsoux p. 94

%

\mbox{}\qed\mbox{}\\

We conclude this section with the definition of a Gelfand triple:
Consider a complete countably normed space $\Phi$, and let
$(\cdot, \cdot)$ denote an inner product on $\Phi$, which is
separately continuous in each variable with respect to the
topology of $\Phi$. Let $\mathcal H$ be the closure of $\Phi$ with
respect to the norm defined by that inner product. The triple
$(\Phi,\mathcal H,\Phi^\prime)$ is called a {\sl Gelfand triple}.
See \cite[p. 101]{MR35:7123}. An important Gelfand triple
consists of the Schwartz space $\sr$ of rapidly decreasing
functions, of the Lebesgue space $\mathbf L_2$ on the real line
and of the Schwartz space of tempered distributions. In the
following section we recall a stochastic counterpart of this
Gelfand triple, which is used below.

\section{The white noise space}

\setcounter{equation}{0} \label{wns} Let $\mathscr S(\mathbb R)$
denote the Schwartz space of {\sl real-valued}, rapidly
decreasing functions. It is a nuclear space, and by the
Bochner-Minlos theorem (see \cite[Th\'eor\`eme 2, p.
342]{MR35:7123}), there exists a probability measure $P$ on the
Borel sets $\mathcal F$ of the dual space $\Omega =\mathscr S(\mathbb
R)^\prime$ such that
\begin{equation}
\label{voltaire} \int_{\Omega}e^{i \langle
\omega ,s\rangle}dP(\omega )= e^{-\frac{\|s\|^2}{2}},\quad \forall
s\in\mathscr S(\mathbb R).
\end{equation}

The real-valued space
\[
\W=
\mathbf L_2(\Omega ,{\mathcal F}, P)
\]
is called the {\sl white noise space}. For $s\in\mathscr S(\mathbb R)$,
let $Q_s$ denote the random variable
\[
Q_s(\omega)=\langle \omega, s\rangle.
\]
It follows from \eqref{voltaire} that
\[
\|s\|_{{\mathbf L}_2(\mathbb R)}=\|Q_s\|_{\W}.\] Therefore, $Q_s$
extends continuously to an isometry from  ${\mathbf L}_2(\mathbb
R)$ into $\W$, which we will still denote by $Q$. In \cite{aal2}
we define
\begin{equation}
\label{opera}
X_m(t)=Q_{T_m(1_{[0,t]})}.
\end{equation}

It follows form the construction in \cite{aal2}, that $X_m(t)$ is
real-valued when $m$ is even. See formula \eqref{eq:Xm} below.\\

In the presentation of the Gelfand triple associated with the
white noise space  we follow \cite{MR1408433}. Let $\ell$ to be
the set of sequences
\begin{equation}
\label{eq:ell}
(\alpha_1,\alpha_2,\ldots),
\end{equation}
indexed by ${\mathbb N}$  with values in ${\mathbb N}_0$, for
which only a finite number of elements $\alpha_j\not=0$. The
white noise space $\W$, being a space of ${\mathbf L}_2$ random
variables on the probability space $(\Omega ,{\mathcal
F}, P)$ specified above, admits a special orthogonal basis
$(H_{\alpha})_{\alpha\in\ell}$,  indexed by the set $\ell$ and
built in terms of the Hermite functions $\widetilde{h_k}$ and of
the Hermite polynomials $h_k$ as
\[
H_\alpha(\omega )=\prod_{k=1}^\infty h_{\alpha_k}(Q_{\widetilde{h_k}}(\omega
)).
\]
We refer the reader to \cite[Definition 2.2.1 p. 19]{MR1408433}
for more information. In terms of this basis, any element $F\in\W$
can be written as
\begin{equation}
\label{fps}
F=\sum_{\alpha\in\ell}f_\alpha H_\alpha,\quad f_\alpha\in\mathbb R,
\end{equation}
with
\[
\|F\|_{\W}^2=\sum_{\alpha\in\W}f_\alpha^2\alpha!<\infty.
\]
There are quite a number of Gelfand triples associated with $\W$.
In \cite{MR2414165}, \cite{al_acap}, and here, we focus on
$(S_1,\W,S_{-1})$, namely the Kondratiev space $S_1$ of
stochastic test functions, $\W$ defined above, and the Kondratiev
space $S_{-1}$ of stochastic distributions. To define these
spaces we first introduce, for $k\in{\mathbb N}$, the Hilbert
space ${\mathcal H}_{k}$ which consists of series of the form
\eqref{fps} such that
\begin{equation}
\label{michelle}
\|F\|_{k}\stackrel{\rm def.}{=}
 \left(\sum_{\alpha\in\ell}(\alpha!)^2f_\alpha^2
(2{\mathbb N})^{k\alpha}\right)^{1/2}<\infty,
\end{equation}
where
\[
(2\mathbb N)^{\pm k\alpha}=(2\cdot1)^{\pm k\alpha_1}(2\cdot 2)^{\pm k\alpha_2}
(2\cdot 3)^{\pm k\alpha_3}\cdots,
\]
and the Hilbert space $\mathcal H^\prime_k$ consisting of
sequences $G=(g_\alpha)_{\alpha\in\ell}$ such that
\[
\|G\|_{k}^\prime\stackrel{\rm def.}{=} \left(\sum_{\alpha\in\ell}g_\alpha^2
(2{\mathbb N})^{-k\alpha}\right)^{1/2}<\infty,
\]
and the duality between an element $F=\sum_{\alpha\in\ell}
f_\alpha H_\alpha\in \mathcal H_k$ and a sequence
$G=(g_\alpha)_{\alpha\in\ell}\in \mathcal H^\prime_k$ is given by
\[
\langle G,F\rangle_{S_{-1},S_1}=\sum_{\alpha\in\ell}\alpha! f_\alpha g_\alpha.
\]

The map which to $F\in\W$ associates its sequence of coefficients
with respect to the basis $(H_\alpha)_{\alpha\in\ell}$ allows to
identify $\W$ as a subspace of $\mathcal H_k^\prime$ for every
$k\in\mathbb N_0$, and it is important to note that
\begin{equation}
\label{ineq}
\|F\|_k^\prime\le\|F\|_{\W},\quad \forall F\in \W.
\end{equation}

The spaces $S_1$ and $S_{-1}$ are defined by
\[
S_1=\bigcap_{k=1}^\infty \mathcal H_k\quad{\rm and}\quad
S_{-1}=\bigcup_{k=1}^\infty \mathcal H^\prime_{k}.
\]

The space $S_1$ is nuclear, see \cite{MR1408433}.\\

The process $\left\{X_m(t)\, ,\,t\in\mathbb R\right\}$ defined in
\eqref{opera}, is written in the series form
\begin{equation}
\label{eq:Xm}
X_m(t)=\sum_{k=1}^\infty\int_0^tT_m\widetilde{h}_k(u)du H_{\epsilon^{(k)}},
\end{equation}
where the series converges in the norm of $\W$, and has an
$S_{-1}$-valued derivative given by the obvious formula
\begin{equation}
\label{Wmseries}
W_m(t)=\sum_{k=1}^\infty (T_m\widetilde{h_k})(t)H_{\epsilon^{(k)}},
\end{equation}
where $\epsilon^{(k)}$ is the sequence in $\ell$ with all entries
equal to $0$, with the exception of the $k$-th, which is equal to
$1$. Furthermore, the series \eqref{Wmseries} converges in the
norm of $\mathcal H_{N+3}^\prime$, where $N$ is as in \eqref{mbound}.
See \cite[Theorem 7.2]{aal2}.

\medskip

{\bf Remark:} Obviously, as $H_{\epsilon^{(k)}}=H_{\epsilon^{(k)}}
(\omega)$, it follows that $X_m(t)=X_m(t,\omega)$, $W_m (t)=W_m
(t,\omega)$. To simplify the notation, we however omit the
$\omega$-dependence throughout, unless specifically required.

\begin{Pn}
We claim:\\

$(a)$  $W_m(t)\in \mathcal H_{N+3}^\prime$ for all $t\in\mathbb
R$.\\

$(b)$ There exists a constant $C_N$ such that
\begin{equation}
\|W_m(t)-W_m(s)\|_{ \mathcal H_{N+3}^\prime}\le C_N|t-s|, \forall
t,s\in \mathbb R. \label{ledru-rollin}
\end{equation}

$(c)$ It holds that
\begin{equation}
\label{deriv}
X_m^\prime(t)=W_m(t),\quad t\in\mathbb R,
\end{equation}
in the norm of $\mathcal H_{N+3}^\prime$, and more generally, in
the norm of any  $\mathcal H_{p}^\prime$ with $p\ge N+3$.
\label{ppn}
\end{Pn}

{\bf Proof:} Claim $(a)$ is proved in \cite[Proof of Theorem 3.2,
p. 1098]{aal2}. It is also shown there, see \cite[Lemma 3.8, p.
1089]{aal2}, that there exist constants $C_1,C_2$ such that
\begin{equation}
\label{manege} \forall t,s\in\mathbb R,\quad
|T_m\widetilde{h}_k(t)-T_m\widetilde{h}_k(s)|\le
|t-s|\cdot(C_1k^{\frac{N+2}{2}}+C_2).
\end{equation}
Since
\[
\|Q_{\widetilde{h_k}}\|_{ \mathcal H_{N+3}^\prime}=(2k)^{-N-3},
\]
we can write $\forall t,s\in\mathbb R$
\[
\begin{split}
\|W_m(t)-W_m(s)\|_{ \mathcal H_{N+3}^\prime}&\le\sum_{k=1}^\infty
|T_m\widetilde{h}_k(t)-T_m\widetilde{h}_k(s)|\|Q_{\widetilde{h_k}}\|_{
\mathcal H_{N+3}^\prime}\\
&\le|t-s|\left\{\sum_{k=1}^\infty(C_1k^{\frac{N+2}{2}}+C_2)(2k)^{-N-3}
\right\}\\
&=C_N|t-s|,
\end{split}
\]
with
\begin{equation}
C_N=\sum_{k=1}^\infty(C_1k^{\frac{N+2}{2}}+C_2)(2k)^{-N-3},
\label{CN}
\end{equation}
which proves $(b)$. We now prove $(c)$. For $t,s\in\mathbb R$,
with $t\not=s$, and $C_N$ as in \eqref{CN}, we have,
\[
\begin{split}
\|\frac{X_m(t)-X_m(s)}{t-s}-W_m(t)\|_{{\mathcal H_{N+3}^\prime}}&=\\
&\hspace{-4cm}=\|\frac{\sum_{k=1}^\infty\int_s^t(
T_m\widetilde{h}_k(u)-T_m\widetilde{h}_k(t))duH_{\epsilon^{(k)}}}{t-s}
\|_{{\mathcal H_{N+3}^\prime}}\\
&\hspace{-4cm}\le  C_N\frac{|\int_s^t|u-t|du}{|t-s|}\\
&\hspace{-4cm}\le  \frac{C_N|t-s|}{2}
\quad{\longrightarrow 0}{{\mbox{ \rm
 as}\,\, s\rightarrow t}}.
\end{split}
\]

The last claim follows from the fact that the spaces $\mathcal
H_n^\prime$ are increasing, with decreasing norms.
\mbox{}\qed\mbox{}\\

The Wick product is defined with respect to the basis
$(H_\alpha)_{\alpha\in\ell}$ by
\[
H_\alpha\lozenge H_\beta=H_{\alpha+\beta}.
\]
It extends to a continuous map from $S_1\times S_1$ into itself
and from $S_{-1}\times S_{-1}$ into itself. Let $l>0$,  and let
$k>l+1$. Consider $h\in {\mathcal H}^\prime_{l}$ and $u\in
{\mathcal H}^\prime_{k}$. Then, V\r{a}ge's inequality holds:
\begin{equation}
\label{vage}
\|h\lozenge u\|_{k}\le A(k-l)\|h\|_{l}\|u\|_{k},
\end{equation}
where
\begin{equation}
\label{vage111}
A(k-l)=\left(\sum_{\alpha\in\ell}(2{\mathbb
N})^{(l-k)\alpha}\right)^{1/2}<\infty.
\end{equation}
See \cite[Proposition 3.3.2, p. 118]{MR1408433}.\\

\bigskip
To conclude this section, we show that the process $X_m$ has $P$-a.s.
continuous sample paths. This property will be utilized in the
last step of the proof of \eqref{oberkampf} below, the Ito formula
associated with $X_m$.

\bigskip

By \cite[Lemma 6.1]{aal2},
\begin{eqnarray}
E [\vert X_m (t) -X_m (s)\vert^2 ]&=&2Re \{ r(t-s)\}\nonumber\\
&\le&2(C_1\vert t-s\vert^2 +C_2 \vert t-s\vert )\nonumber\\
&\le& 2(C_1 +C_2 )(\vert t-s\vert^2 \vee\vert t-s\vert ),
\label{E7.0}
\end{eqnarray}
where the first equality is due to item (2) of aforementioned Lemma,
while the following inequality is due to item (3), with $C_1$, $C_2$
some finite, positive constants.

Recall that $X_m$ is a Gaussian process. Then, for all $t,s\in{\mathbb R}$,
$\vert t-s\vert \le 1$, it follows from (\ref{E7.0}) that
\begin{equation}
E[ \vert X_m (t) -X_m (s) \vert^4 ]\le 12 (C_1 +C_2 )^2 \vert t-s\vert^2.
\end{equation}
By Kolmogorov's continuity criterion, see e.g. \cite[Theorem I-1.8]{Revuz91},
it follows that there exists a $P$-a.s. continuous modification of
$X_m$, a modification we consider here.

\section{The Wick-Ito integral}

\setcounter{equation}{0}
\label{wickito}
The main result of this section is the following theorem:

\begin{Tm}
\label{bastille}
Let $Y(t)$, $t\in[a,b]$ be an $S_{-1}$-valued
function, continuous in the strong topology of $S_{-1}$. Then,
there exists a $p\in\mathbb N$ such that the function $t\mapsto
Y(t)\lozenge W_m(t)$ is $\mathcal H_p^\prime$-valued, and
\begin{equation*}
\int_a^b Y(t,\omega)\lozenge W_m(t)dt
=\lim_{\left|\Delta\right|\to 0} \sum_{k=0}^{n-1}
Y(t_k,\omega)\lozenge \left(X_m(t_{k+1})-X_m(t_k)\right),
\end{equation*}
where the limit is in the $\mathcal H_p^\prime$ norm, with
$\Delta: a=t_0<t_1<\cdots <t_n=b$ a partition of the interval
$[a,b]$ and $\left|\Delta\right|=\max_{0\leq k \leq
n-1}(t_{k+1}-t_k)$.
\end{Tm}

{\bf Proof:} We proceed in a number of steps.\\

STEP 1: {\sl $W_m(t)\in\mathcal{H}_{N+3}^\prime$ for $t\in\mathbb
R$, and  satisfies \eqref{ledru-rollin}:
\[
\|W_m(t)-W_m(s)\|_{{\mathcal H}^\prime_{N+3}}\le C_N|t-s| ,\quad \forall
t,s\in\mathbb R.
\]
for some constant $C_N$.}\\

See Proposition \ref{ppn}.\\

STEP 2:  {\sl There exists a $p\in\N$, $p>N+3$, such that
$Y(t)\in\mathcal{H}_{p}^\prime$ for all $t\in [a,b]$, being
uniformly continuous from $[a,b]$ into
$\mathcal{H}_{p}^\prime$.}\\

Theorem \ref{function} with $E=[a,b]$ ensures that a $p$ (not
necessarily larger than $N+3$) exists with the stated properties.
Since the norms $\|\cdot\|_{{\mathcal H_p^\prime}}$ are
decreasing, we may assume  that $p>N+3$.

%But if $Y$ is $\mathcal H_p$-valued and continuous in

%the $\mathcal H_p$-norm, the same hold for any $r>p$ since the

%norms are decreasing.

Using V\r{a}ge's inequality \eqref{vage}, it follows that, for
$p>N+3$,
\[
\begin{split}
\left\|Y(t)\lozenge W_m(t)-Y(s)\lozenge W_m(s)\right\|_p&\le\\
&\hspace{-4cm}\le \left\|(Y(t)-Y(s))\lozenge W_m(t)\right\|_p+
\left\|Y(s)\lozenge (W_m(t)-W_m(s))\right\|_p\\
&\hspace{-4cm}\le A(p-N-3)\|Y(t)-Y(s)\|_p\|W_m(s)\|_{N+3}+\\
&\hspace{-3.5cm}+A(p-N-3)\|Y(s)\|_p\|W_m(t)-W_m(s)\|_{N+3}
%\leq A(q) \left\|Y\right\|_p \left\|W_m\right\|_3 ??
%\left\|Y\lozenge W_m\right\|_q^2\leq A(q) \left\|Y\right\|_p
%\left\|W_m\right\|_3
\end{split}
\]
where $A(p-N-3)$ is defined by \eqref{vage111}, with
$\|\cdot\|_p\stackrel{\rm def.}{=} \|\cdot\|_{{\mathcal
H}_p^\prime}$ used to simplify the notation.\\

In view of Step 2, the integral $\int_a^b Y(t)\lozenge W_m(t)dt$
makes sense as a Riemann integral of a Hilbert space valued
continuous function.\\

STEP 3: {\sl Let $\Delta$ be a partition of the interval $[a,b]$.
We now compute an estimate for}
\[
\begin{split}
\int_a^b Y(t)\lozenge W_m(t)dt-\sum_{k=0}^{n-1} Y(t_k)\lozenge
\left(X_m(t_{k+1})-X_m(t_k)\right)&=\\
&\hspace{-6cm}= \sum_{k=0}^{n-1}
\left(\int_{t_k}^{t_{k+1}}(Y(t)-Y(t_k)) \lozenge W_m(t)dt\right).
\end{split}
\]

Let $p$ be as in Step 2, and let $\epsilon>0$. Since $Y$ is
uniformly continuous on $[a,b]$ there exists an $\eta>0$ such that
\[
|t-s|<\eta\Longrightarrow \|Y(t)-Y(s)\|_p<\epsilon.
\]

Set
\[
\tilde{C}=\max_{s\in[a,b]}\|W_m(s)\|_{N+3}\quad{\rm and}\quad
A=A(p-N-3).
\]

Let $\Delta$ be a partition of $[a,b]$ with
\[
|\Delta|=\max\left\{|t_{k+1}-t_k|\right\}<\eta.
\]

We then have:

\[
\begin{split}
\left\|\sum_{k=0}^{n-1} \left(\int_{t_k}^{t_{k+1}}(Y(t)-Y(t_k))
\lozenge W_m(t)dt \right)\right\|_p\leq\\
&\hspace{-7cm}\le \sum_{k=0}^{n-1}
\left(\int_{t_k}^{t_{k+1}}\left\|(Y(t)-Y(t_k)) \lozenge
W_m(t)\right\|_p dt\right)\\
&\hspace{-7cm}\leq A\sum_{k=0}^{n-1}
\left(\int_{t_k}^{t_{k+1}}\left\|(Y(t)-Y(t_k))\right\|_p\left\|
W_m(t)\right\|_{N+3} dt\right)\\
&\hspace{-7cm}\leq \tilde{C}A\sum_{k=0}^{n-1}\int_{t_k}^{t_{k+1}}
\left\|(Y(t)-Y(t_k))\right\|_pdt\\
&\hspace{-7cm}\le \epsilon\tilde{C}A(b-a)
\end{split}
\]
\mbox{}\qed\mbox{}\\

\section{An Ito formula}

\setcounter{equation}{0}
\label{itoform}
We extend the classical Ito's formula to the present setting. We
need the extra assumption that the function
\[
r(t)=\|T_m1_{[0,t]}\|_{{\mathbf L_2(\mathbb R)}}
\]
is absolutely continuous with respect to the Lebesgue measure.
This is in particular the case for the fractional Brownian motion.
This is also the case e.g. for the function $m$ defined in
\eqref{u4}, where, for that $m$,
\[
r(t)=
\frac{\sqrt{2\pi}}{8}\left\{
1-e^{-\frac{t^2}{8}}(1+t^2)\right\}.
\]

\begin{Tm}
Suppose that $r(t)$ is absolutely continuous with respect to the
Lebesgue measure. Let $f:\R\to\R$ be a $C^2(\mathbb R)$ function.
Then
%\[
%f(X_m(t)),\quad \int_0^t f^\prime(X_m(s))\lozenge W_m(s)ds,\quad
%\int_0^tf^{\prime\prime} (X_m(s)) r(s)^\prime ds
%\]
%all belong to ${\mathcal H}^\prime_{N+3}$ for every $t\in[c,d]$.
%Furthermore,
\begin{equation}
\label{oberkampf}
\begin{split}
f(X_m(t))&=f(X_m (t_0 ))+\int_{t_0}^t f'(X_m(s))\lozenge W_m(s)ds+\\
&\hspace{5mm}+ \frac{1}{2}\int_{t_0}^tf^{\prime\prime}(X_m(s))
r^\prime(s)
ds,\quad t_0 <t \in\R,
\end{split}
\end{equation}
where the equality is in the $P$-almost sure sense.
%${\mathcal H}^\prime_{N+3}$.
\end{Tm}

{\bf Proof:}  We prove for $t>t_0 =0$. The proof for any other
 interval in $\R$ is essentially the same.
We divide the proof into a number of steps. Step
1-Step 8 are constructed so as to show that \eqref{oberkampf}
holds, $\forall t >0$, for  $C^2$ functions with
compact support, with the equality holding in the ${\mathcal
H}^\prime_p$ sense. This implies its validity in the $P$-a.s.
sense (actually, holding  $\forall \omega\in\Omega$),  hence,
setting the ground for the concluding Step 9, in which the result
is extended
to hold for all $C^2$  functions $f$.\\

STEP 1: {\sl For every $(u,t)\in\mathbb R^2$, it holds that
\[
e^{iuX_m(t)}\in\mathcal W,
\]
and
\begin{equation}
e^{iuX_m(t)}\lozenge W_m(t)\in\mathcal H^\prime_{N+3}.
\label{leipzig2010}
\end{equation}

}

Indeed, since $X_m$ is real, we have
\[
|e^{iuX_m(t)}|\le 1,\quad \forall u,t\in\mathbb R,
\]
and hence $e^{iuX_m(t)}\in\W$. Since $\W\subset \mathcal
H_{N+1}^\prime$ and since $W_m(t)\in \mathcal H_{N+3}^\prime$ for
all $t\in\mathbb R$, it follows from  V\r{a}ge's inequality
\eqref{vage} that \eqref{leipzig2010} holds.\\

In the following two steps, we prove formula \eqref{oberkampf} for
exponential functions. For $\alpha\in\R$ we set:
\[
g(x)=\exp(i\alpha x).
\]

STEP 2: {\sl It holds that
\begin{equation}
\label{gprimem}
 g'(X_m(t))=i\alpha g(X_m(t))\lozenge
W_m(t)+\frac{1}{2} (i\alpha)^2g(X_m(t))r^\prime (t).
\end{equation}}

Indeed, $g(X_{m}(t))$ belongs to $S_{-1}$, see \cite[p.
65]{MR1408433}, and we have from \cite[Lemma 2.6.16, p.
66]{MR1408433}:
\[
\begin{split}
g(X_m(t))=\exp(i\alpha X_m(t))&=\exp^\lozenge
\left(i\alpha X_m(t)
+\frac{1}{2}(i\alpha)^2
\left\|T_mI_t\right\|_{\mathbf L_2(\R)}^2\right)\\
&=\exp^\lozenge \left(i\alpha X_m(t) +\frac{1}{2}(i\alpha)^2 r(t)\right).
\end{split}
\]

%STEP 2: {\sl The function $t\mapsto Y\lozenge W_m$ is continuous

%on $[a,b]$}\\

The hypothesis that $r$ is absolutely continuous with respect to
Lebesgue measure comes now into play. Since the function
$t\mapsto W_m(t)$ is continuous in $S_{-1}$ and since
$X_m^\prime=W_m$, see Proposition \ref{ppn},  an application of
\cite[Theorem 3.1.2]{MR1408433} with
\[
h(t)=i\alpha W_m(t)-\frac{\alpha^2}{2}r^\prime (t)
\]
leads to
\[
\begin{split}
g'(X_m(t))&
%=\exp^\lozenge(i\alpha X_m(t)
%+\frac{1}{2}(i\alpha)^2\textsl{Re}~r(t)) \lozenge(i\alpha W_m(t)
%+\frac{1}{2}(i\alpha)^2 \left({\rm Re}~r(t)
%\right)')\\
%&=
=g(X_m(t))\lozenge(i\alpha W_m(t)+\frac{1}{2}(i\alpha)^2
r(t)^\prime) \\
&= g(X_m(t))\lozenge(i\alpha W_m(t))+\frac{1}{2}
(i\alpha)^2g(X_m(t))r^\prime(t).
\end{split}
\]
We thus obtain \eqref{gprimem}.\\

STEP 3: {\sl Equation \eqref{oberkampf} holds for exponentials.}\\

Indeed, it follows from  \eqref{gprimem} that:
\[
\begin{split}
g(X_m(t))&=g(X_m(0))+\int_0^ti\alpha g(X_m(s))\lozenge
W_m(s)ds\\
&\hspace{5mm}+\frac{1}{2}\int_0^t(i\alpha)^2g(X_m(s))
r^\prime (s) ds.
\end{split}
\]
This can be written
\[
\begin{split} g(X_m(t))&=g(0)+\int_0^tg'(X_m(s))\lozenge
W_m(s)ds+\\
&\hspace{5mm}+\frac{1}{2}\int_0^tg''(X_m(s)) r^\prime (s) ds,
\end{split}
\]
that is
\begin{equation}
\label{montmartre}
\begin{split}
e^{iuX_m(t)}=1&+\int_0^tiue^{iuX_m(s)}\lozenge W_m(s)ds\\
&+\frac{1}{2}\int_0^t(iu)^2e^{iuX_m(s)} r^\prime (s)ds
\end{split}
\end{equation}

In the following two steps, we prove \eqref{oberkampf} to hold
for Schwartz functions.\\

STEP 4: {\sl The function $(u,t)\mapsto e^{iuX_m(t)}\lozenge
W_m(t)$ is continuous from $\mathbb R^2$ into $\mathcal
H^\prime_{N+3}$.}\\

We first recall that the function $t\mapsto X_m(t)$ is
continuous, and even uniformly continuous, from $\mathbb R$ into
$\mathcal W$, and hence from $\mathbb R$ into $\mathcal
H_{N+5}^\prime$ since
\[
\|u\|_{\mathcal H_{N+3}^\prime}\le\|u\|_{\W},\quad{\rm for}\quad
u\in\W.
\]
Therefore the function $(u,t)\mapsto e^{iuX_m(t)}$ is continuous
from $\mathbb R^2$ into $\mathcal H^\prime_{N+3}$. Furthermore,
\[
\begin{split}
\|e^{iu_1X_m(t_1)}\lozenge W_m(t_1)-e^{iu_2X_m(t_2)}\lozenge
W_m(t_2)\|_{\mathcal H^\prime_{N+3}}&\le\\
&\hspace{-25mm}\le \| (e^{iu_1X_m(t_1)}-e^{iu_2X_m(t_2)})
\lozenge W_m(t_1)\|_{\mathcal H^\prime_{N+1}}+\\
&\hspace{-20mm}+
 \| e^{iu_1X_m(t_1)}\lozenge(W_m(t_2)-
W_m(t_1))\|_{\mathcal H^\prime_{N+3}}\\
&\hspace{-65mm}\le A(2)\| (e^{iu_1X_m(t_1)} -e^{iu_2X_m(t_2)})
\|_{\mathcal H^\prime_{N+1}}\cdot\| W_m(t_1)\|_{\mathcal H^\prime_{N+3}}+\\
&\hspace{-60mm}+ A(2) \| e^{iu_1X_m(t_1)}\|_{\mathcal
H^\prime_{N+1}}\cdot\|(W_m(t_2)- W_m(t_1))\|_{\mathcal
H^\prime_{N+3}},
\end{split}
\]
where $A(2)$ is defined by \eqref{vage111}. This completes the
proof of STEP 4 since $t\mapsto W_m(t)$ is continuous in the norm
of $\mathcal H_{N+3}^\prime$ and $(u,t)\mapsto e^{iuX_m(t)}$ is
continuous in the norm of
$\mathcal H_{N+1}^\prime$.\\

STEP 5: {\sl  \eqref{oberkampf} holds for $f$ in the Schwartz
space.}\\

%\[

%\begin{split}

%f(X_m(t))=1&+\int_0^t\left(\int_\R iu h(u)

%e^{iuX_m(s)}du\right)\lozenge W_m(s)ds\\

%&+\frac{1}{2}\int_0^t\left(\int_\R(iu)^2h(u)e^{iuX_m(s)}du\right)

%r(s)^\prime ds.

%\end{split}

%\]

%}

Let $s$ be in the Schwartz space. Replace $u$ by $-u$ in
\eqref{montmartre},  and multiply both sides of this equation by
$s(u)$. Integrating with respect to $u$, and interchanging order
of integration, we obtain
\[
\begin{split}
\int_{\mathbb R} s(u)e^{-iuX_m(t)}du&=\int_{\mathbb R} s(u)du+\\
&\hspace{5mm}+ \int_0^t\left(\int_{\mathbb R}(-iu)
s(u)e^{-iuX_m(s)}du\right)
\lozenge W_m(s)ds+\\
&\hspace{5mm}+\frac{1}{2}\int_0^t\int_{\mathbb
R}\left(\int_0^t(-iu)^2 s(u)du\right)e^{-iuX_m(s)} r(s)^\prime ds.
\end{split}
\]
Continuity of the function in the previous step allowed to use
Fubini's theorem for functions with values in a Hilbert space
(see \cite[Theorem 2.6.14, p.65]{DS1}, \cite[Proposition 9, p.
97]{bourbaki_int_5}), and to interchange the
order of integration.

%Since $X_m$ is real valued,  we can write:

%\[

%f(X_m(t))=\int_a^b  h(u)e^{iuX_m(t)}du,

%\]

%and this function is bounded in absolute value, and therefore

%belongs to ${\mathbf L}_2(dP)$.

Since
\[
\widehat{s}^{\,\,\prime}(x)=\int_{\mathbb
R}(-iu)s(u)e^{-iux}du,\quad{\rm and}\quad
\widehat{s}^{\,\,\prime\prime}(x)=\int_{\mathbb
R}(-iu)^2s(u)e^{-iux}du,
\]
we obtain
\[
\begin{split}
\widehat{s}(X_m(t))=&\widehat{s}(0)+\int_0^t
(\widehat{s})^\prime(X_m(s))\lozenge W_m(s)ds\\
&+\frac{1}{2}\int_0^t (\widehat{s})^{\prime\prime}(X_m(s))
r(s)^\prime ds.
\end{split}
\]
This completes the proof of STEP 5 since the Fourier transform maps the Schwartz
space onto itself.\\

To show that \eqref{oberkampf} holds for $f$ of class $C^2$ with
compact support, we will use the concept of approximate identity.
For $\epsilon>0$, define
\[
k_\epsilon(x)=\frac{1}{\sqrt{2\pi}\epsilon}\exp-{\frac{x^2}{2\epsilon^2}}.
\]

STEP 6: {\sl It holds that
\begin{equation}
\int_{\mathbb R} k_\epsilon(x)dx=1,
\label{int1}
\end{equation}
and, for every $r>0$
\begin{equation}
\label{int2}
\lim_{\epsilon\rightarrow
0}\int_{|x|>r}k_\epsilon(x)dx=0.
\end{equation}}

Indeed, \eqref{int1} follows directly from the fact that $k_\epsilon$ is an
$\mathcal N(0,\epsilon^2)$ density.
%\[

%\int_{\mathbb R}e^{-x^2}dx=\sqrt{2\pi}.

%\]

Furthermore, for $|x|>r>0$,
\[
\begin{split}
\frac{1}{\epsilon\sqrt{2\pi}}\int_{r}^\infty
e^{-\frac{x^2}{2\epsilon^2}}dx&=
\frac{1}{\epsilon\sqrt{2\pi}}\int_{r}^\infty \frac{x}{\epsilon^2}
e^{-\frac{x^2}{2\epsilon^2}}\frac{\epsilon^2}{x}dx\\
&\le \frac{\epsilon}{r\sqrt{2\pi}} \int_{r}^\infty
\frac{x}{\epsilon^2} e^{-\frac{x^2}{2\epsilon^2}}dx\\
&=\frac{\epsilon}{r\sqrt{2\pi}}e^{-\frac{r^2}{2\epsilon^2}}\\
&\longrightarrow 0\quad\mbox{\rm as}\quad \epsilon\rightarrow 0.
\end{split}
\]

The properties in STEP 6 express the fact that $k_\epsilon$ is an
approximate identity. Therefore, it follows from \cite[Theorem
1.2.19, p. 25]{MR2463316} that, for every continuous function
with compact support,
\[
\lim_{\epsilon\rightarrow 0}\|k_\epsilon \ast f-f\|_\infty=0.
\]

STEP 7: {\sl The functions
\[
(k_\epsilon \ast f)(x)=\frac{1}{\sqrt{2\pi}\epsilon}\int_{\mathbb
R} \exp (-{\frac{(u-x)^2}{\epsilon^2}} )f(u)du
\]
are in the Schwartz space.}\\

One proves by induction on $n$ that the $n$-th derivative
\[
(k_\epsilon \ast f)^{(n)}(x)
\]
is a finite sum of terms of the form
\[
\frac{1}{\sqrt{2\pi}\epsilon}\int_{\mathbb R}
\exp \big( -{\frac{(u-x)^2}{\epsilon^2}}\big)p(x-u)f(u)du,
\]
where $p$ is a polynomial. That all limits,
\[
\lim_{|x|\rightarrow\infty}x^m
(k_\epsilon \ast f)^{(n)}(x)=0
\]
is then shown using the dominated convergence theorem.\\

STEP 8: {\sl \eqref{oberkampf} holds for $f$ of class $C^2$ and
with compact
support.}\\

A function $f$ of class $C^2$ with compact support can
be approximated, together with first two derivatives,  in the
supremum norm by Schwartz functions. This is done as follows.
Take for simplicity $\epsilon=\frac{1}{n}$, $n=1,2,\ldots$. We apply
\cite[Theorem 1.2.19, p. 25]{MR2463316} to $f$,
$f^\prime$ and $f^{\prime\prime}$. Set
\[
a_n=k_{1/n}\ast f,\quad  b_n=k_{1/n}\ast f^\prime,\quad{\rm
and}\quad c_n=k_{1/n}\ast f^{\prime\prime}.
\]

Integration by parts shows that
\[
\begin{split}
a_n^\prime&=b_n\\
b_n^\prime&=c_n.
\end{split}
\]

Furthermore,
\[
\begin{split}
\lim_{n\rightarrow\infty}\|a_n-f\|_\infty&=0,\\
\lim_{n\rightarrow\infty}\|b_n-f^\prime\|_\infty&=0,\\
\lim_{n\rightarrow\infty}\|c_n-f^{\prime\prime}\|_\infty&=0.
\end{split}
\]

For every $n$ we have:
\[
a_n(X_m(t))=a_n(0)+\int_0^ta_n^\prime(X_m(s))\lozenge
W_m(s)ds+\frac{1}{2}\int_0^ta_n^{\prime\prime}(X_m(s))r^\prime(s)ds.
\]

We claim that, for a given $t$, the sequence
$(a_n(X_m(t)))_{n\in\mathbb N}$ is a Cauchy sequence in any
$\mathcal H^{\prime}_p$ since
\[
\|a_n(X_m(t))-a_m(X_m(t))\|_{\mathcal
H^{\prime}_p}\le\|a_n(X_m(t))-a_n(X_m(s))\|_{\W}\le\|a_n-a_m\|_\infty,
\]
and denote by $f(X_m(t))$ the corresponding limit.\\

Similarly, the sequence
\[
(\int_0^t b_n(X_m(u))\lozenge W_m(u)du)_{n\in\mathbb N}
\]
is a Cauchy sequence in $\mathcal H^{\prime}_p$, since
\[
\begin{split}
\|\int_0^tb_n(X_m(u))\lozenge W_m(u)du-\int_0^tb_m(X_m(u))\lozenge
W_m(u)du\|_{\mathcal
H^{\prime}_p}&\le\\
&\hspace{-4cm}\le\int_0^t\|(b_n-b_m)(X_m(u))\lozenge W_m(u)\|_{\mathcal
H^{\prime}_p}du\\
&\hspace{-4cm}\le A(2)\int_0^t\|(b_n-b_m)(X_m(u))\|_{\mathcal
H^{\prime}_p} \|W_m(u)\|_{\mathcal H^{\prime}_{N+3}}
du\\
&\hspace{-4cm}\le
A(2)\|b_n-b_m\|_\infty\times\int_0^t\|W_m(u)\|_{\mathcal
H^{\prime}_{N+3}}du.
\end{split}
\]

Denote
\[
\int_0^tf^\prime(X_m(u))\lozenge W_m(u)du
\]
to be its limit. A similar argument holds for
\[
\int_0^t c_n(X_m(u))r^\prime(u)du.
\]
Details are omitted.\\

Note that we have actually shown \eqref{oberkampf} to hold for
$C^2$ functions with compact support with the equality understood
in the $\mathcal H_p^\prime$ sense. This implies that it also
holds in the $P$-a.s sense.\\

STEP 9: {\sl \eqref{oberkampf} holds with probability $1$ for all
$f\in C^2(\mathbb R)$.}\\

Here, we follow key arguments
of the corresponding proofs of the standard Ito rule,
given in e.g., \cite[Theorem
IV-3.3, p. 138] {Revuz91},
\cite[Theorem 3.3, p. 149]{MR917065}. Specifically, the following
standard localization argument is utilized.
Let  $\tau_N$ be a stopping time defined by
\[
\tau_N=\inf\left\{s>0\,:\, |X_m(s)|>N\right\}.
\]

Set
\[
X_m^N(s)\stackrel{\rm def.}{=}X_m(s\wedge\tau_N).
\]

Then, by Step 8, \eqref{oberkampf} holds for $ \left\{X_m^N(s)\,
,\, s\ge 0\right\}$, a.s., for all $f\in C^2(\mathbb R)$. Fix an
arbitrary $\epsilon>0$ and let $N=N(\epsilon)<\infty$ be such that
\[
P(\sup_{0\le s\le t}|X_m(s)|>N)<\epsilon.
\]
As \eqref{oberkampf} holds for $X_m^N$ a.s., it then follows that
\eqref{oberkampf} holds for
$X_m$ with probability greater that $1-\epsilon$, for all $f\in
C^2(\mathbb R)$.  The arbitrariness of  $\epsilon$ completes the proof
of the fact that  \eqref{oberkampf} holds, $P$-a.s., for all $t\in \mathbb R$. This suggests that both sides of \eqref{oberkampf}
are {\it modifications} of one another. Since both are $t$-continuous
(see the discussion at the end of Section 4 for the continuity of the LHS), they are in fact
{\it indistinguishable} processes, which is to say that \eqref{oberkampf} holds for all
$t\in \mathbb R$,
$P$-a.s.
\mbox{}\qed\mbox{}\\

\section{Concluding Remarks}

\setcounter{equation}{0} {\bf 1.} Note that no adaptability of
the integrand with respect to an  underlying filtration is
assumed. In this sense, one may regard the integral defined here
in fact as a Wick-Skorohod
integral.\\

{\bf 2.} Due to the fact that  $\| F\|_p^{\prime} \le\| F\|_{\W}$
for all $F\in\W$ (see \eqref{ineq}),   it follows that the
integral defined in Theorem  \ref{bastille}, being an ${\mathcal
H}_p^{\prime}$ limit of Riemann sums, is defined in a weaker
sense than the standard  Ito and Skorohod integrals, which are
defined in an $\mathbf L_2$ sense. This is a reasonable price to
pay so as to integrate with respect to a larger class of
non-$\mathbf L_2$ integrands. This places the proposed
integral well within existing stochastic integration theory,
as a non-trivial extension is formulated, at the (expected)
expense of a somewhat weaker sense of convergence.\\

{\bf  3.} We note the reduction of the calculus derived here to the
standard Ito calculus when $r(t)=|t|$. This case corresponds to
setting $H=1/2$ in \eqref{BH} and \eqref{VH}, so that $V_{1/2}=1$
and $m(u)=\frac{1}{2\pi}$. For example, for $Y(t)=B(t)$, both stochastic
integrals give
\[
\int_0^tB(t)dB(t)=\frac{B^2(t)-t}{2}.
\]
See \cite{MR2414165}. Furthermore, for the fractional Brownian
motion with Hurst parameter $H\in(0,1)$, i.e. (up to a
multiplicative constant) $r(t)=|t|^{2H}$, our integral coincides
with that proposed in e.g. \cite{MR1741154} and \cite{MR1956473},
\cite{MR2046814}.\\

{\bf 4.} Note that our Ito formula for $X_m$ being the fractional
Brownian motion with Hurst parameter $H\in(0,1)$, hence
$m(t)=\frac{1}{2\pi}|t|^{1-2H}$, coincides with that of Bender
\cite{MR1956473} specified for $C^2$ functions with values in a
space of distributions. We note however the difference between
the proofs. In \cite{MR1956473}
%our Ito formula and that
%derived by Bender for the fractional Brownian motion
%To extend the Ito formula to functions which do not have
%necessarily a compact support, one can proceed as in
%\cite{MR1956473}. There,
Bender shows that the $S$-transforms of both sides of equation
\eqref{oberkampf} agree. The conclusion that
\eqref{oberkampf} holds in fact for all $\omega\in\Omega$ follows from
the fact that the $S$-transform is injective. Here one can use
Bender's approach for $C^2$ distributions (with derivatives
understood in the sense of distributions), replacing Bender's
$r^\prime(t)=t^{2H-1}$ ($t>0$) with the derivative of a general
$r$. We omit the computations which are essentially the same.
% The

%utilization of the $S$-transform is however restricted to

%distribution functions

\bibliographystyle{plain}
\def\cprime{$'$} \def\lfhook#1{\setbox0=\hbox{#1}{\ooalign{\hidewidth
  \lower1.5ex\hbox{'}\hidewidth\crcr\unhbox0}}} \def\cprime{$'$}
  \def\cprime{$'$} \def\cprime{$'$} \def\cprime{$'$} \def\cprime{$'$}

\end{document}